\documentclass[a4paper,10pt]{article}
\usepackage[utf8]{inputenc}

\usepackage{mathtools}
\usepackage{amssymb}
\usepackage{amsthm}
\usepackage{amsmath}
\usepackage[bottom=4.5cm]{geometry}
\usepackage{authblk}
\usepackage{hyperref}

\newtheorem{theorem}{Theorem}
\newtheorem{corollary}{Corollary}
\newtheorem{lemma}{Lemma}
\newtheorem{definition}{Definition}
\newcommand{\F}{\mathbb{F}}

\renewcommand\footnotemark{}
\title{Generalized Almost Perfect Nonlinear Binomials and Trinomials Over Fields of Prime-Square Order\thanks{This manuscript version has been accepted for publication, after peer review
but is not the Version of Record and does not reflect post-acceptance improvements, or any
corrections. The Version of Record is available online at: \url{https://doi.org/10.1016/j.ffa.2023.102185}.}}
\author{Christof Beierle}
\affil{Faculty of Computer Science, Ruhr University Bochum, Bochum, Germany\\ \href{mailto:christof.beierle@rub.de}{christof.beierle@rub.de}}
\date{}

\begin{document}
\maketitle

\begin{abstract}
Let $p>3$ be a prime. We show that, for each integer $d$ with $p \leq d \leq 2(p-1)$, there exists a generalized almost perfect nonlinear (GAPN) binomial or trinomial over $\F_{p^2}$ of algebraic degree $d$. We start by deriving sufficient conditions for the function $G \colon \F_{p^2} \rightarrow \F_{p^2}, X \mapsto X^{d_1} + u X^{d_2}$ to be GAPN in the case where one of the terms of $G$ is GAPN. We then give explicit constructions of GAPN binomials over $\F_{p^2}$ of any odd algebraic degree between $p$ and $2(p-1)$ and, in the case where $p$ is not a Mersenne prime, also of any even algebraic degree in this range. To obtain GAPN functions of even algebraic degree also in the general case, we finally show how to construct GAPN trinomials over $\F_{p^2}$ of any even algebraic degree between $p$ and $2(p-1)$ by applying a characterization of a special form of GAPN binomials by \"{O}zbudak and S\u{a}l\u{a}gean. Our constructed functions are the first GAPN functions of even algebraic degree over extension fields of odd characteristic reported so far.

\end{abstract}

\section{Introduction}
In the following, let $p$ be an odd prime and let $n$ be a positive integer. By $\F_{p^n}$ we denote the finite field with $p^n$ elements, and by $\F_{p^n}^*$ we denote its multiplicative group. We study \emph{generalized almost perfect nonlinear (GAPN) functions}, which were introduced by Kuroda and Tsujie as a generalization of APN functions in even characteristic.
\begin{definition}[Kuroda and Tsujie, \cite{DBLP:journals/ffa/KurodaT17}]
A function $G \colon \F_{p^n} \rightarrow \F_{p^n}$ is \emph{GAPN} if, for all elements $a \in \F_{p^n}^*$ and $b \in \F_{p^n}$, the equation 
\[ D_a^{(p-1)}G(X) \coloneqq \sum_{i \in \F_p} G(X + ia) = b\]
has at most $p$ solutions $X \in \F_{p^n}$.
\end{definition}

As it was outlined in~\cite{DBLP:journals/ffa/OzbudakS21}, the object $D_a^{(p-1)}G$ equals the $p-1$-th order discrete derivative of $G$ at point $a$, and we simply call $D_a^{(p-1)}G$ a \emph{derivative} in the following. Thus, for a function to be GAPN, any derivative at $a \neq 0$ has to be $p$-to-1 (see~\cite[Prop.\@ 1]{DBLP:journals/ffa/OzbudakS21}). For $p=2$, this definition coincides with the one for APN functions, see, e.g.,~\cite[Def.\@ 2.5]{DBLP:journals/dcc/Pott16}. While APN functions, especially APN permutations, have applications in cryptography (they have in fact been originally introduced by Nyberg and Knudsen in~\cite{DBLP:conf/crypto/NybergK92} as permutations providing optimal resistance against differential cryptanalysis), it was shown that GAPN functions (resp., APN functions) have applications in finite geometry as some of them can be used to construct dual arcs~\cite{DBLP:journals/ffa/KurodaT17} (resp., dual hyperovals~\cite{yoshiara2008dimensional}). Note that the case of $n=1$ is trivial, i.e., \emph{any} function over $\F_p$ is GAPN. The reason is that any derivative is constant, and thus, a $p$-to-1 function. We will therefore always assume $n \geq 2$ in the following. 

An important notion for functions over finite fields is the notion of the \emph{algebraic degree}. Let $G \colon \F_{p^n} \rightarrow \F_{p^n}, X \mapsto \sum_{u=0}^{p^n-1} a_u X^u$ with $a_u \in \F_{p^n}$. For a non-negative integer $u$, we define $d_p(u) \coloneqq \sum_{i=0}^{\infty} h_i$, where $\sum_{i=0}^{\infty} h_i p^i$ is the base-$p$ representation of $u$. If $G \neq 0$, the \emph{algebraic degree} of $G$ is defined as the maximum value of $ d_p(u)$ over all $u \in \{0,\dots,p^n-1\}$ with $a_u \neq 0$.

In~\cite{DBLP:journals/ffa/OzbudakS21}, the authors introduced the notion of generalized EA-equivalence, which is an invariant for the GAPN property. Two functions $F,G \colon \F_{p^n} \rightarrow \F_{p^n}$ are \emph{generalized EA-equivalent} if there exist affine permutations $A,B$ of $\F_{p^n}$ and a function $C \colon \F_{p^n} \rightarrow \F_{p^n}$ of algebraic degree at most $p-1$ such that $G = B \circ F \circ A + C$. Up to this notion of equivalence, the only GAPN functions over $\F_{p^n}, n \geq 2$, known so far are monomial functions (\cite{DBLP:journals/ijfcs/Kuroda20,DBLP:journals/ffa/KurodaT17,DBLP:journals/ffa/OzbudakS21,ffa_new_monomial,DBLP:journals/ffa/ZhaHZ18}), i.e., functions of the form $X \mapsto X^d$, linear combinations of GAPN monomials (\cite{DBLP:journals/ffa/KurodaT17,DBLP:journals/ffa/OzbudakS21}), linear combinations of monomials of the form $X^{tp^{n-1}-1}$ with $t$ being even (\cite[Prop.\@ 4.3]{DBLP:journals/ffa/ZhaHZ18}), or functions of algebraic degree $p$ (see \cite{DBLP:journals/ffa/KurodaT17,DBLP:journals/ffa/OzbudakS21}). This raises the question on the existence of GAPN functions that are inequivalent to functions coming from those classes, and in particular on the possible algebraic degrees of GAPN functions. Note that the algebraic degree, if it is at least $p$, is invariant under generalized EA-equivalence. For $n \geq 2$, we know that $p$ is a lower bound on the algebraic degree of GAPN functions (simply because for function of algebraic degree strictly less than $p$, any derivative is constant) and that $n(p-1)$ is a trivial upper bound (simply because a function over $\F_{p^n}$ cannot have higher algebraic degree).

The simplest type of GAPN functions are GAPN monomial functions and several results on them have been reported in the literature: 
Kuroda proved that GAPN monomials over finite extension fields of odd characteristic cannot have even algebraic degree~\cite{DBLP:journals/ijfcs/Kuroda20}, so there are no GAPN functions of even algebraic degree over extension fields reported in the literature yet (besides the APN case of $p=2$). For $n \geq 2$, the maximum algebraic degree of GAPN monomials is therefore $n(p-1)-1$ and it is attained by the permutation $X \mapsto X^{p^n-2}$, i.e., the permutation fixing zero and mapping each element in $\F_{p^n}^*$ to its multiplicative inverse. Moreover, Kuroda proved that any monomial function of algebraic degree $n(p-1)-1$ is GAPN~\cite{DBLP:journals/ijfcs/Kuroda20}. In~\cite{DBLP:journals/ffa/BartoliGPZ22}, the authors associated monomial functions to algebraic curves and were able to study GAPN monomials by methods from algebraic geometry.  From~\cite{DBLP:journals/ffa/OzbudakS21},  we know the following for GAPN monomials with exponents of the form $kp^{r_1} + \ell p^{r_2}$ with $k,\ell \in \{0,\dots,p-1\}$ and $r_1,r_2 \in \{0,\dots,n-1\}$.

\begin{lemma}[\"{O}zbudak and S\u{a}l\u{a}gean, \footnote{In a recent talk at the \emph{7th International Workshop on
Boolean Functions and their Applications (BFA) 2022}, \"{O}zbudak and S\u{a}l\u{a}gean  presented a complete classification of GAPN monomials with exponents of the form $kp^{r_1} + \ell p^{r_2}$.  For $n>2$, they showed that the (sufficient) Condition 1 is also necessary and, for $n=2$, they showed that the (necessary) Condition 2 is also sufficient.} 
\cite{DBLP:journals/ffa/OzbudakS21}]
\label{lem:special_gapn_monomial}
Let $n$ be a positive integer and let $d = k p^{r_1} + \ell p^{r_2}$ with $k,\ell \in \{0,\dots,p-1\}$, $r_1,r_2 \in \{0,\dots,n-1\}$, $r_1 \neq r_2$, and $p \leq k+\ell < 2(p-1)$.
\begin{enumerate}
\item If $\gcd(r_1-r_2,n) = 1$ and $\gcd(k + \ell - (p-1),p^n-1) = 1$, then $X \mapsto X^d$ is a GAPN function over $\F_{p^n}$.
\item If $X \mapsto X^d$ is a GAPN function over $\F_{p^n}$, we have $\gcd(r_1-r_2,n)=1$ and $\gcd(k+\ell-(p-1),p^{n_1}-1)=1$ for all positive integers $n_1$ with $n_1 \mid n$ and $n_1 < n$.
\end{enumerate}
\end{lemma}

Note that in the case of a finite field of prime-square order, i.e., $n=2$, all exponents of monomial functions must be of the form $kp + \ell$ with $k,\ell \in \{0,\dots,p-1\}$. In particular, besides the fact that we cannot obtain GAPN monomials of even algebraic degree, Lemma~\ref{lem:special_gapn_monomial} implies that we also cannot always obtain GAPN monomials of \emph{every odd} algebraic degree between $p$ and $2(p-1)$. For instance, if $n=2, p=11$, we cannot reach algebraic degree 15, simply because $\gcd(15-(p-1),p-1) = \gcd(5,10) = 5 \neq 1$.

In this note, we focus on the case of $n=2$ and first study binomial functions $X \mapsto X^{d_1} + u X^{d_2}$, where one of the terms is a GAPN monomial. In Section~\ref{sec:main_construction}, we provide some sufficient conditions when such a binomial is GAPN and in Section~\ref{sec:explicit_families}, we provide explicit constructions of GAPN binomials of \emph{any} algebraic degree $d$ with $p \leq d \leq 2(p-1)$ in the case where $p$ is not a Mersenne prime.\footnote{See Sequence \href{https://oeis.org/A000668}{A000668} in the On-Line Encyclopedia of Integer Sequences (OEIS) for a list of the first Mersenne primes.} If $p$ is a Mersenne prime, we can still construct GAPN binomials of any \emph{odd} algebraic degree in this range. In Section~\ref{sec:trinom}, we apply a characterization by \"{O}zbudak and S\u{a}l\u{a}gean of a special form of GAPN binomials consisting of GAPN terms of the same algebraic degree in order show a sufficient condition for the trinomial function $X \mapsto u \cdot X^{2p-1} + v \cdot X^{3p-2} + X^{h(p+1)}$ to be GAPN. For $p>3$, this then allows us to show the existence of GAPN trinomials of every even algebraic degree between $p$ and $2(p-1)$ without any further restriction on $p$.

\section{Constructing GAPN Binomials with One GAPN Term}
\label{sec:main_construction}
In this section, we prove the following theorem, which allows to construct GAPN binomials over fields of prime-square order from GAPN monomials. We denote a monomial function $X \mapsto X^d$ shortly by $M_d$.

\begin{theorem}
\label{thm:main}
Let $p$ be an odd prime and $d_1,d_2 \in \{1,\dots,p^2-1\}$. Let $M_{d_1}$ be a GAPN function over $\F_{p^2}$ and let $u \in \F_{p^2}^*$.  The function
\[G \colon \F_{p^2} \rightarrow \F_{p^2}, \quad X \mapsto X^{d_1} + u \cdot X^{d_2}\ \]
is GAPN if one of the following two conditions hold:
\begin{enumerate}
    \item $d_2$ is odd and $u$ is not a square in $\F_{p^2}$.
    \item $d_2$ is even and there exists an odd integer $N$ such that $p + 1$ is a multiple of $N$, $d_2 -d_1$ is a multiple of $N$, and $u$ is not an $N$-th power in $\F_{p^2}$.
\end{enumerate}
\end{theorem}

For the case of monomials over $\F_{p^2}$, we have the following identity, which will be at the core of our proofs.
\begin{lemma}\label{lem:power}
Let $d \in \{0,1,\dots,p^2-1\}$. Over $\F_{p^2}$, we then have
\[\left(D_1^{(p-1)} M_d(X)\right)^p = (-1)^{d}D_1^{(p-1)} M_d(X). \]
\end{lemma}
\begin{proof}
Let us fix an element $x \in \F_{p^2}$. Using the fact that the Frobenius mapping $X \mapsto X^p$ is linear and fixes each element of $\F_p$, we obtain
\begin{align*}
    \left(D_1^{(p-1)} M_d(x)\right)^p = \left( \sum_{i \in \F_p} (x + i)^d\right)^p = \sum_{i \in \F_p}(x^p + i)^d.
\end{align*}
Now, since we are in the case of extension degree 2 over $\F_p$, we have $x^p+x = tr(x) \in \F_p$, where $tr$ denotes the absolute trace function. Therefore, we have
\[ \sum_{i \in \F_p}(x^p + i)^d = \sum_{i \in \F_p} (-x + tr(x) + i)^d = (-1)^d\sum_{i \in \F_p} (x - tr(x) - i)^d = (-1)^d D_1^{(p-1)}M_d(x).\]
\end{proof}

 For any $a \in \F_{p^n}^*$, we know (see~\cite{DBLP:journals/ffa/ZhaHZ18}) that $D_a^{(p-1)} M_d(X) = a^d D_1^{(p-1)}M_d(a^{-1}X)$, which yields
\begin{align*}
    \left( D_a^{(p-1)}M_d(X)\right)^p &= a^{dp} \left( D_1^{(p-1)}M_d(a^{-1}X)\right)^p = (-1)^{d}a^{dp}D_1^{(p-1)}M_d(a^{-1}X) \\&= (-1)^{d}A^{d}D_a^{(p-1)}M_d(X),
\end{align*}
where $A \coloneqq a^{p-1}$. We now have all the ingredients to prove Theorem~\ref{thm:main}.

\begin{proof}[Proof of Theorem~\ref{thm:main}]
Let us consider a binomial $G \colon X \mapsto X^{d_1} + u \cdot X^{d_2}$ over $\F_{p^2}$, where $1 \leq d_1,d_2 \leq p^2-1, u \in \F_{p^2}^*$ and let us assume that $M_{d_1}$ is a GAPN monomial. Let $a \in \F_{p^2}^*$. For $G$ to be GAPN, we need to show that the derivative $D_a^{(p-1)}G$ is $p$-to-1. We have
\begin{equation}\label{eq:D_G} D_a^{(p-1)}G(X) =  D_a^{(p-1)} M_{d_1}(X) + u \cdot D_a^{(p-1)} M_{d_2}(X)
\end{equation}
and raising it to the $p$-th power yields
\begin{equation}\label{eq:D_G_p} \left(D_a^{(p-1)}G(X)\right)^p =  (-1)^{d_1}A^{d_1}D_a^{(p-1)} M_{d_1}(X) + u^p \cdot (-1)^{d_2} A^{d_2} D_a^{(p-1)} M_{d_2}(X),
\end{equation}
where $A \coloneqq a^{p-1}$.
Multiplying Equation~(\ref{eq:D_G}) by $(-1)^{d_2}u^{p-1}A^{d_2}$ and subtracting Equation~(\ref{eq:D_G_p}) yields
\[ \left((-1)^{d_2}u^{p-1}A^{d_2}\right) D_a^{(p-1)}G(X) - \left(D_a^{(p-1)}G(X)\right)^p = H \cdot D_a^{(p-1)} M_{d_1}(X),\]
where $H \coloneqq (-1)^{d_2}u^{p-1}A^{d_2} - (-1)^{d_1}A^{d_1}$. If $H \neq 0$, the mapping $H \cdot D_a^{(p-1)} M_{d_1}$ is $p$-to-$1$, because $M_{d_1}$ is a GAPN monomial by assumption. In that case, it immediately follows that $D_a^{(p-1)}G$ is $p$-to-$1$, i.e., $G$ is GAPN. Indeed, assume on the contrary that there exist pairwise distinct $x_1,\dots,x_{p+1} \in \F_{p^2}$ such that, for all $1 \leq i,j \leq p+1$, $D_a^{(p-1)}G(x_i)=D_a^{(p-1)}G(x_j)$, we would get $H \cdot D_a^{(p-1)} M_{d_1}(x_i) = H \cdot D_a^{(p-1)} M_{d_1}(x_j)$ for all $1 \leq i,j \leq p+1$, a contradiction to the $p$-to-$1$ property of $H \cdot D_a^{(p-1)} M_{d_1}$. To derive conditions on $G$ being GAPN, it therefore suffices to derive the conditions on $H$ being nonzero.  

Since $M_{d_1}$ is GAPN, we know that $d_1$ must be odd, hence $(-1)^{d_1} = -1$.
We then have
\begin{align*}
    H &= -A^{d_1} \left((-1)^{d_2+1} u^{p-1}A^{d_2 -d_1} - 1 \right),
\end{align*}
which is zero if and only if \begin{equation}\label{eq:main_condition}(-1)^{d_2 + 1} (u a^{d_2-d_1})^{p-1} = 1.\end{equation}

\paragraph{Case $d_2$ is odd.}
In this case, Equation~(\ref{eq:main_condition}) is equivalent to $u a^{d_2-d_1} \in \F_p$. Since $d_1$ is odd, this implies $u a^{2m} \in \F_p$ for some $m$. But then, $u$ must be a square in $\F_{p^2}$ since all elements in $\F_p$ are squares in $\F_{p^2}$ and, obviously, $a^{2m}$ is a square.   Thus, if $u$ is not a square in $\F_{p^2}$, the function $G$ is GAPN.

\paragraph{Case $d_2$ is even.}
Equation~(\ref{eq:main_condition}) is equivalent to $u^{p-1} a^{(p-1)(d_2-d_1)} = (-1)$. Suppose that $d_2-d_1$ is a multiple of $N$ for an odd integer $N$, this implies $u^{p-1}$ to be an $N$-th power since $(-1)$ is an $N$-th power (viz., $(-1)^N = -1$). We have $u^{p-1} = u^{-2} u^{p+1}$, so if $p = -1 \mod N$, the element $u^2$ must be an $N$-th power. But then, also $u$ must be an $N$-th power. Indeed, if $u = g^e$ for a generator $g \in \F_{p^2}^*$ and $2e = 0 \mod N$, we have $e = 0 \mod N$ since $N$ is odd. Therefore, $G$ is GAPN if there exists an odd integer $N$ such that $p = -1 \mod N$, $d_2-d_1 = 0 \mod N$, and $u$ is not an $N$-th power. 
\end{proof}

We stress that the conditions given in Theorem~\ref{thm:main} are not necessary to obtain GAPN binomials. For instance, the function $X \mapsto X^{3\cdot 7+4} + X^{6 \cdot 7 + 4}$ over $\F_{7^2}$ is GAPN. Although $X \mapsto X^{3 \cdot 7 + 4}$ is a GAPN monomial, this binomial cannot be from the construction of Theorem~\ref{thm:main}. Indeed, for $a \in \F_{7^2}$ with minimal polynomial $X^2+1$, we have $a^{(7-1)(3 \cdot 7)} = (-1)$, which implies Equation~(\ref{eq:main_condition}) (and thus $H=0$) to hold for that particular value of $a$. Moreover, there exist GAPN binomials where none of its terms is GAPN itself. For example, if $u$ is a primitive element of $\F_{11^2}^*$, the function $X \mapsto X^{2 \cdot 11 + 10} + u \cdot X^{5 \cdot 11 + 10}$ is GAPN over $\F_{11^2}$.

\subsection{Some Explicit Families}
\label{sec:explicit_families}
The most simple GAPN monomial is the generalized Gold function $X \mapsto X^{2p-1}$ (see~\cite{DBLP:journals/ffa/KurodaT17}), which is of algebraic degree $p$. We base our first two constructions on this monomial.
Our goal is to construct GAPN functions of every algebraic degree $d$ with $p \leq d \leq 2(p-1)$. The case of odd algebraic degree is straightforward. 
\begin{corollary}
\label{cor:odd}
Let $p$ be an odd prime, $u \in \F_{p^2}^*$ be a primitive element and let $k,\ell \in \{0,\dots,p-1\}$ with $k+\ell$ being odd. Then, the function
\[G_{k,\ell} \colon \F_{p^2} \rightarrow \F_{p^2}, X \mapsto X^{2p-1} + u \cdot X^{kp+ \ell}\]
is GAPN and of algebraic degree $\max \{p, k+\ell\}$.
\end{corollary}

It is more complicated to derive constructions for the case of even algebraic degree, and we only manage to construct GAPN binomials of even algebraic degree in the case where $p$ is not a Mersenne prime, i.e., $p \neq 2^n-1$ for $n$ being a natural number. Let us first give a simple construction in the case where $p = 2 \mod 3$.

\begin{corollary}
Let $p$ be an odd prime with $p = 2 \mod 3$ and let $u$ be a primitive element of $\F_{p^2}^*$. Let $h \in \{1,\dots, p-1\}$. Then,
\[ G_h \colon \F_{p^2} \rightarrow \F_{p^2}, X \mapsto X^{2p-1} + u \cdot X^{h(p+1)}\]
is GAPN and of algebraic degree $\max \{p,2h\}$.
\end{corollary}
\begin{proof}
We have $p = -1 \mod 3$ and $h(p+1) - (2p-1) = (h-2)p + h + 1 = 0 \mod 3$. By a similar argument as before, a primitive element of $\F_{p^2}^*$ cannot be a cube in $\F_{p^2}$. The GAPN property of $G_h$ then follows from Theorem~\ref{thm:main}.
\end{proof}

A slightly more complicated construction can be given in the general case, as shown below. If $p$ is odd and not a Mersenne prime, we can thus construct GAPN binomials over $\F_{p^2}$ of any algebraic degree between $p$ and $2(p-1)$.
\begin{corollary}
\label{thm:even_construction}
Let $p$ be an odd prime which is not Mersenne. Then, for any even integer $d$ with $p \leq d \leq 2(p-1)$, there exists a GAPN binomial over $\F_{p^2}$ of algebraic degree $d$.
\end{corollary}
\begin{proof}
Since $p$ is not a Mersenne prime, we have $p + 1 = m N$ for a positive integer $m$ and an odd integer $N \geq 3$ and thus, $p = -1 \mod N$. Let $u$ be a primitive element of $\F_{p^2}^*$. Since $p^2-1 = (p-1)(p+1)$, we have $\gcd(N,p^2-1) = N > 1$, so the mapping $X \mapsto X^N$ is not a permutation of $\F_{p^2}^*$, which implies that $u$ is not an $N$-th power. Let now $h \in \{1,\dots,p-1\}$ and let $d_1 \coloneqq k_1 p + \ell_1, d_2 \coloneqq h(p+1)$ with $k_1 = \frac{p-N}{2}, \ell_1 = \frac{p+N}{2}$. Then, $d_2$ is even and 
\[d_2 - d_1 = \left(\frac{p-N}{2} - \frac{p+N}{2}\right) \mod N= 0,  \]
so $d_2 -d_1$ is a multiple of $N$. Further, we have 
\[ k_1 + \ell_1 - (p-1) = p -(p-1) = 1,\]
so $\gcd(k_1+\ell_1-(p-1),p^2-1) = 1$ and $X \mapsto X^{k_1p + \ell_1}$ is GAPN over $\F_{p^2}$ by Lemma~\ref{lem:special_gapn_monomial}. By Theorem~\ref{thm:main}, the function $G \colon X \mapsto X^{d_1} + u \cdot X^{d_2}$ is GAPN over $\F_{p^2}$ and it is of algebraic degree $\max\{k_1 + \ell_1, 2h \} = \max \{ p,2h\}$. Putting in the values for $d_1,d_2$, the function $G$ can be expressed as
\[G \colon \F_{p^2} \rightarrow \F_{p^2}, X \mapsto X^{\frac{p^2+p -N(p-1)}{2}} + u \cdot X^{h(p+1)}.\]
\end{proof}

Given a (possibly large) prime $p$, the element $N$ for the construction in Corollary~\ref{thm:even_construction} can be efficiently computed with a complexity of $\mathcal{O}(\log_2 p)$ integer divisions by repeatedly dividing $p+1$ by 2, until the result is odd.

The question on the possible algebraic degrees of GAPN binomials over $\F_{p^2}$ with $p$ being a Mersenne prime remains open. Certainly, there exist such primes for which we cannot have GAPN binomials of all even algebraic degrees between $p$ and $2(p-1)$. For $p=3$, a computer search yields that there are no GAPN functions over $\F_{p^2}$ of even algebraic degree. For $p=7$, our search yields that there do not exist GAPN binomials over $\F_{p^2}$ of algebraic degree 8 or 12, even if we drop the condition that one of the terms is GAPN. But there exist GAPN binomials over $\F_{7^2}$ of algebraic degree 10, e.g., the example mentioned after the proof of Theorem~\ref{thm:main}. 

\section{GAPN Trinomials of Even Algebraic Degree}
\label{sec:trinom}
To show the existence of GAPN functions of every even algebraic degree between $p$ and $2(p-1)$ also in the case where $p > 3$ is a Mersenne prime, we will construct GAPN trinomials where two of its terms are GAPN and of algebraic degree $p$ and the remaining term is of even algebraic degree. To do so, we use the following characterization of a special form of GAPN polynomials.

\begin{lemma}[\"{O}zbudak and S\u{a}l\u{a}gean, see Thm.\@ 5 in \cite{DBLP:journals/ffa/OzbudakS21}]
\label{lem:binoms}
Let $n$ be a positive integer, $s \in \{1,\dots,p-2\}$ with $\gcd(s,p^n-1)=1$, and let $r_1,r_2 \in \{0,\dots,n-1\}, r_1 < r_2$ with $\gcd(r_2-r_1,n) = 1$. For $i \in \{s,s+1,\dots,p-1\}$, let $d_i \coloneqq i p^{r_2} + (p-1+s-i)p^{r_1}$. Let $k$ be the inverse of $\frac{p^{k_2}-p^{k_1}}{p-1} \mod p^n-1$.

For $a \in \F_{p^n}^*$, the derivative $D_{a^k}^{(p-1)}G$ of the function $G \colon X \mapsto \sum_{i=s}^{p-1} c_i X^{d_i}$, $c_i \in \F_{p^n}$, is $p$-to-1 if and only if 
\begin{equation*}
\label{eq:binom_derivative_condition}
    \sum_{i=0}^{p-1-s}c_{i+s} \binom{p-1-s}{i}(-a^{p-1})^i \neq 0.
\end{equation*}
\end{lemma}

In our case, we have $n=2, r_1=0, r_2=1$, and therefore $k=1$. To simplify our construction, we will apply Lemma~\ref{lem:binoms} only to the case of $s=1$ and $c_i = 0$ for $i \notin \{1,2\}$. We obtain that, for $a \in \F_{p^2}^*$, the derivative $D_a^{p-1}G$ of the function $G \colon X \mapsto c_1X^{p + (p-1)} + c_2 X^{2p + (p-2)}$ is $p$-to-1 if and only if 
\begin{equation}
\label{eq:binom_derivative_condition_special}
    c_1 + 2c_2a^{p-1} \neq 0.
\end{equation}

Using Equation~(\ref{eq:binom_derivative_condition_special}), we are now able to construct GAPN trinomials as follows. Fixing a primitive element $g$ of $\F_{p^2}^*$, we denote by $\langle g^{p-1} \rangle$ the set of $(p-1)$-th powers in $\F_{p^2}^*$.
\begin{theorem}
\label{thm:trinom}
Let $p$ be an odd prime and let $u, v \in \F_{p^2}^*$ such that $2vX^5 + uX^4 + u^pX + 2v^p \in \F_{p^2}[X]$ has no root in $\langle g^{p-1} \rangle$. Then, for any $h \in \{0,\dots,p-1\}$, the function 
\[G_h \colon \F_{p^2} \rightarrow \F_{p^2}, \quad X \mapsto u \cdot X^{2p-1} + v \cdot X^{3p-2} + X^{h(p+1)}\]
is GAPN and of algebraic degree $\max\{p,2h\}$.
\end{theorem}
\begin{proof}
Let $d_1 = 2p-1 = p + (p-1)$, $d_2 = 3p-2 = 2p + (p-2)$, and $d_3 = h(p+1)$. Let $a \in \F_{p^2}^*$ and define $A \coloneqq a^{p-1}$. For $G_h$ to be GAPN, we need to show that the derivative $D_a^{(p-1)}G_h$ is $p$-to-1. We have 
\begin{equation}\label{eq:trinomial_derivative_1}D_a^{(p-1)}G_h(X) = u \cdot D_a^{(p-1)}M_{d_1}(X) + v \cdot D_a^{(p-1)}M_{d_2}(X) + D_a^{(p-1)} M_{d_3}(X) \end{equation}
and, by Lemma~\ref{lem:power}, raising it to the $p$-th power yields
\begin{equation}\label{eq:trinomial_derivative_2}\left(D_a^{(p-1)}G_h(X)\right)^p = -u^p \cdot A^{d_1} D_a^{(p-1)}M_{d_1}(X) - v^p \cdot A^{d_2} D_a^{(p-1)}M_{d_2}(X) + A^{d_3} D_a^{(p-1)} M_{d_3}(X). \end{equation}
Since $A^{d_3} = a^{(p^2-1)h} = 1$, by subtracting Equation~(\ref{eq:trinomial_derivative_2}) from Equation~(\ref{eq:trinomial_derivative_1}) we obtain
\[D_a^{(p-1)}G_h(X) - \left(D_a^{(p-1)}G_h(X)\right)^p = (u+u^pA^{d_1})D_a^{(p-1)}M_{d_1}(X) + (v + v^pA^{d_2})D_a^{(p-1)}M_{d_2}(X),\]
and, similarly as in the proof of Theorem~\ref{thm:main}, it suffices to show that the mapping $H \coloneqq (u+u^pA^{d_1})D_a^{(p-1)}M_{d_1} + (v + v^pA^{d_2})D_a^{(p-1)}M_{d_2}$ is $p$-to-1. Indeed, $H$ is the derivative at point $a$ of the binomial function $X \mapsto (u+u^pA^{d_1})X^{d_1} + (v + v^pA^{d_2})X^{d_2}$, which is (according to Equation~(\ref{eq:binom_derivative_condition})) $p$-to-1 if and only if 
\begin{equation}\label{eq:final}u + u^p A^{d_1} + 2 vA + 2v^pA^{d_2+1} \neq 0.\end{equation}
We have $A^{d_1} = A^{p+(p-1)} =A^{2(p+1)}A^{-3}=A^{-3}$ and $A^{d_2+1} = A^{2p + (p-2)+1} = A^{3(p+1)}A^{-4} = A^{-4}$, hence, Equation~(\ref{eq:final}) is equivalent to 
\[ uA^4 + u^pA + 2vA^5 + 2v^p \neq 0.\]
The $p$-to-1 property of $H$ follows since the polynomial $2vX^5 + uX^4 + u^pX + 2v^p \in \F_{p^2}[X]$ has no root in $\langle g^{p-1} \rangle$ by assumption. 
\end{proof}

The condition on the polynomial $2vX^5 + uX^4 + u^pX + 2v^p$ given in Theorem~\ref{thm:trinom} is independent on $h$. Therefore, to show the existence of GAPN functions of \emph{every} even algebraic degree between $p$ and $2(p-1)$, it is enough to show the existence of a polynomial of the above form with no roots in $\langle g^{p-1} \rangle$. We do so in the next lemma (setting $v = 2^{-1} \in \F_p^*$).
\begin{lemma}
Let $p>3$ be a prime. There exists an element $u \in \F_{p^2}^*$ such that $X^5 + u X^4 + u^pX + 1 \in \F_{p^2}[X]$ has no root in $\langle g^{p-1} \rangle$.
\end{lemma}
\begin{proof}
In the following, we will assume $p \geq 31$. The case of $p \in \{ 5,7,11,13,17,19,23,29\}$ was handled by a computer search. For $X \in \F_{p^2}$, let $S_X \coloneqq \{ u \in \F_{p^2}^* \mid  X^5 + uX^4 + u^pX + 1 = 0\}$ and let $S \coloneqq \bigcup_{X \in \langle g^{p-1} \rangle} S_X$. We will show that $\lvert S \rvert  < p^2-1$.

We first observe that, for each element $X \in \langle g^{p-1} \rangle$, the set $S_X$ contains at most $p$ elements since it consists of the (non-zero) roots of a polynomial of degree $p$. If $X = -1$, we have $X^5 + uX^4 + u^pX + 1 = u-u^p$, which is zero if and only if $u \in \F_{p}$. Hence, $S_{-1} = \F_{p}^*$. Moreover, we can deduce that, for each element $X \in \langle g^{p-1} \rangle$ with $X^3+1 \neq 0$, the set $S_X$ contains an element $ \mu_X  \in \F_p^*$, namely $\mu_X = -\frac{X^5+1}{X^4+X}$. Indeed, for $X \in \langle g^{p-1} \rangle$, we have $X^p = X^{-1}$, hence
\begin{align*}
    (X^5+1)(X^4+X)^p &= X^{4p+5} + X^{p+5} + X^{4p} + X^p = X + X^4 + X^{5p+1} + X^{5p+4} \\ &= (X^4+X)(X^5+1)^p
\end{align*}
and it follows that $\mu_X = \mu_X^p$.

Further, for $X \in \langle g^{p-1}\rangle$, the element $\nu_X \coloneqq -X^{-4}$ is a member of both $S_X$ and $S_{-X}$, which implies $\lvert S_X \cup S_{-X} \rvert \leq 2p-1$. If $X^8-1 \neq 0$, the element $\nu_X$ is not in $\F_p$. Hence, if $X^8-1 \neq 0, X^3+1 \neq 0$ and $X^3-1 \neq 0$, we have $\lvert (S_X \cup S_{-X}) \setminus \F_p^* \rvert \leq 2p-3$. In all other cases, we have $\lvert (S_X \cup S_{-X}) \setminus \F_p^* \rvert \leq 2p-1$. Since the solutions of $X^8 = 1$ come in pairs $(X,-X)$, there are at most three sets $\{X,-X\}$ with $X \in \langle g^{p-1} \rangle, X \neq \pm 1$, $X^8-1 = 0$. Additionally, there are at most two sets $\{X,-X\}$ with $X \in \langle g^{p-1} \rangle, X \neq \pm 1$, $X^3+1 = 0$ and at most two sets $\{X,-X\}$ with $X \in \langle g^{p-1} \rangle, X \neq \pm 1$, $X^3-1 = 0$. This yields at most 7 such exceptional pairs among all $X \in \langle g^{p-1}\rangle \setminus \{-1,1\}$.  

For a primitive element $g \in \F_{p^2}^*$, we have $\langle g^{p-1}\rangle = \{ g^{i(p-1)} \mid 0 \leq i <p+1\} = \{\pm g^{i(p-1)} \mid 0 \leq i \leq \frac{p-1}{2}\}$. We finally obtain
\begin{align*} \lvert S \rvert &\leq \lvert S_{-1} \rvert + \lvert S_{1} \setminus \F_p^*\rvert + \sum_{i=1}^{\frac{p-1}{2}} \lvert (S_{g^{i(p-1)}} \cup S_{-g^{i(p-1)}})\setminus \F_p^* \rvert \\ 
&\leq (p-1) + (p-1) + \frac{p-1}{2}(2p-3) + 7 \cdot 2 =p^2 - \frac{p}{2} - \frac{1}{2} + 14 < p^2-1,\end{align*}
where the last inequality holds because $p \geq 31$. 
\end{proof}

\section{Concluding Remarks}
Our technique of constructing GAPN binomials  is not applicable for higher extension degrees of $\F_p$. For example, in the proof of Theorem~\ref{thm:main}, we cancel one of the terms $t_1(X), t_2(X)$ of the derivative in Equation~(\ref{eq:D_G}) by taking its $p$-th power and then expressing each of the resulting terms $t_1(X)^p$, resp., $t_2(X)^p$ by $c_1 t_1(X)$, resp., $c_2 t_2(X)$, where $c_1,c_2$ are some constants. A similar approach is used in the proof of Theorem~\ref{thm:trinom}. Unfortunately, for $n>2$, we have the following fact.
\begin{lemma}
Let $n > 2$ and let $G \colon \F_{p^n} \rightarrow \F_{p^n}$ be a function for which there exists $r \in \{1,\dots,n-1\}$ and an element $c \in \F_{p^n}$ such that $(D_1^{(p-1)}G(X))^{p^r} = c \cdot D_1^{(p-1)}G(X)$ for all $X \in \F_{p^n}$. Then, $G$ is not GAPN.
\end{lemma}
\begin{proof}
Let us fix an element $c \in \F_{p^n}$ such that $(D_1^{(p-1)}G(X))^{p^r} = cD_1^{(p-1)}G(X)$ holds for all $X \in \F_{p^n}$.
Since $0 < r < n$, the equation $Y^{p^r} - cY = 0$ has at most $p^r$ solutions $Y \in \F_{p^n}$. By taking $Y \coloneqq D_1^{(p-1)}G(X)$, this implies $\lvert \mathrm{Im}(D_1^{(p-1)}G)\rvert \leq p^r$. Thus, $D_1^{(p-1)}G$ can only be $p$-to-1 if $r=n-1$. But then, we have $Y^{p^{n-1}} - cY = 0$ if and only if $Y - c^pY^p = 0$, which has at most $p$ solutions $Y \in \F_{p^n}$. Hence, $\lvert \mathrm{Im}(D_1^{(p-1)}G)\rvert \leq p$. So, $D_1^{(p-1)}G$ can only be $p$-to-1 if $n=2$.
\end{proof}

The question on the possible algebraic degrees of GAPN functions over $\F_{p^n}, n>2$ remains open. Note that even in the special case of $p=2$ this question is open (see the discussion in~\cite{DBLP:journals/tit/BudaghyanCHLS18}) and probably very hard to settle completely. Moreover, for extensions of $\F_p$ of degree higher than two and with $p$ being odd, it is left as an open problem to construct GAPN functions of even algebraic degree.

\paragraph{Acknowledgment.}
The author thanks the reviewers for their detailed and useful comments, which helped improving the presentation of the results.

The author is  funded  by Deutsche  Forschungsgemeinschaft  (DFG) under Germany's Excellence Strategy - EXC 2092 CASA - 390781972.

\end{document}